# Group Invariant Bounded Linear Functions on Dedekind Complete Totally Ordered Riesz Spaces

George Chailos

*Department of Mathematics, University of Nicosia, 1700, Nicosia, Cyprus*
chailos.g@unic.ac.cy, telephone no: 0035799354236

**Abstract**

In this paper we consider the set *B* of all bounded subsets of *V*, where *V* is a totally ordered Dedekind complete Riesz space equipped with the order topology. We show the existence of *bounded linear functions* on *B* that are invariant under group actions of the symmetric group of $B$. To do this, we construct a set of "approximately" group invariant bounded linear functions and we show, using Tychonff's Theorem (that is equivalent to the Axiom of Choice), that this set has a cluster point. This cluster point is the group invariant bounded linear function on *B* that we are looking for.

**1. Introduction**

In this paper we extend a result proven in [2]. There, the authors proved the existence of shift invariant mean-like functions on the set *B* of all bounded subsets of totally ordered Dedekind complete Riesz spaces. Here we show the existence of *bounded linear functions* that are invariant under group actions of the symmetric group of $B$. We also note that in [2] the proof required only a special case of Tychonff's theorem which is equivalent to the Ultrafilter Theorem[1]. The same holds for here, but for matters of simplicity and of not repeating the same arguments as in [2] we employ the full version of the Tychonff's theorem (which is equivalent to the Axiom of Choice, see [7]).

A *Riesz space V* (see [6]) is a vector space endowed with a partial order "$\preceq$" that

(a) If $x \preceq y$, then $x + z \preceq y + z$, $\forall x, y, z \in V$

(b) If $x \preceq y$, then $\alpha x \preceq \alpha y$, $\forall x, y \in V$ and $\forall 0 \leq \alpha \in \mathbb{R}$.

(c) For all $x, y \in V$, their "supremum" $x \vee y$ and their "infimum" $x \wedge y$ with respect to "$\preceq$" exist and are elements of $V$.

---

[1] The authors of [2] showed that this special case of Tychonoff's theorem (which holds for Hausdorff compact spaces) is a consequence of the Ultrafilter theorem, which is strictly weaker than the Axiom of Choice [4].

Here we let $V$ to be a totally ordered Dedekind complete Riesz space. The total ordering ensures that all elements are comparable with respect to $\preceq$ and hence the order topology is well defined on $V$. The Dedekind completeness means that every subset of $V$ which is bounded above has a least upper bound (supremum). We set the order topology $T_{ord}$ on $V$. Thus, $(V, T_{ord})$ is in a sense a generalization of $(\mathbb{R}, T_{ord})$. Hence, $(V, T_{ord})$ becomes a regular space (we refer to [8] for details regarding the order topology). It is well known (see [1]) that every element $x$ in a Riesz space $V$ has unique positive $x^+ = x \vee 0$ and negative $x^- = -x \vee 0$ parts. Also, $x = x^+ - x^-$ and an absolute value is defined as $|x| = x^+ + x^-$. Regarding convergence in a Riesz space $V$, a net $\{x_\alpha\}$ in $V$ is said to converge monotonely, if it is a monotone decreasing (respectively increasing) net and its infimum (respectively supremum) $x$ exists in $V$ and is denoted by $x_\alpha \downarrow x$ (respectively $x_\alpha \uparrow x$). A net $\{x_\alpha\}$ in a Riesz space $V$ is said to *converge in order* to $x$ if there exists a net $\{a_\alpha\}$ in $V$ such that $a_\alpha \downarrow 0$ and $|x_n - x| \preceq a_\alpha$.

**Example 1.1:** Consider $\mathbb{R}^n$ with the usual vector operations, its usual norm derived by the order topology, and the partial order defined by:
$$x \preceq y \Leftrightarrow x_k \leq_\mathbb{R} y_k, k = 1, 2, ..., n.$$
Then, $\mathbb{R}^n$ is a totally ordered Dedekind complete Riesz space.

It is worth noting that we can also define a norm $\|\cdot\|$ on a Riesz space $V$ in such a way that it becomes a *normed lattice* $(V, \|\cdot\|)$. That is, if $x \preceq y$, then $\|x\| \leq \|y\|$, $\forall x, y \in V$.
A normed lattice which is also a Banach space is called a *Banach Lattice* (see [1]).

**Example 1.2:** All of the classical Banach Spaces, $l_p$, $c$, $c_o$, $C(K)$, $L^p(\mu)$ are Banach lattices in their usual norm and the pointwise (almost everywhere in the case of $L^p(\mu)$) order.

The next result shows that there exists a relation between $(V, T_{ord})$ and $(V, \|\cdot\|)$. The proof of it and hence the result holds for order intervals in general partially ordered Riesz spaces.

**Lemma 1.3**
The closed order intervals of $(V, T_{ord})$ are also closed in $(V, \|\cdot\|)$.

*Proof:* At first we show that the positive cone $(V^+, T_{ord})$ is normed closed. To see this, assume that a net $\{y_\lambda\} \in V^+$ satisfies $y_\lambda \to y$ in $(V, \|\cdot\|)$. We have to show that $y \in V^+$. Indeed, since $|y_\lambda - y^+| \preceq |y_\lambda - y| \in V^+$ and since $\|\cdot\|$ is a lattice norm, $\|y_\lambda - y^+\| \leq \|y_\lambda - y\|$ and hence $y_\lambda \to y^+$.

From the uniqueness of the limit point (since $V, \|\cdot\|$ is Hausdorff), $y = y^+$. Hence, $y \in V^+$. Thus, $V^+$ is closed. Now let $x, y \in V$ and assume w.l.o.g. that $x \preceq y$. Then, it is easy to observe (see [6]) that the order closed interval $[x, y]$ is written as $[x, y] = (x + V^+) \cap (y - V^+)$. Since $(V^+, T_{ord})$ is normed closed, so are the sets $(x + V^+), (y - V^+)$, and hence their intersection. Thus, $[x, y]$ is normed closed. □

## 2. Group invariant Bounded Linear Functions

The question that we answer affirmatively here is the following: Is there any natural way of associating a *Bounded Linear Function (B.L.F.)* to the set $B$ of bounded subsets of $(V, T_{ord})$ which is invariant under *group actions on B*?

**Formulation of the problem**: Consider doubly infinite sequences $\{x_n\}_{n \in \mathbb{Z}}$ in $(V, T_{ord})$, that is, elements of $V^{\mathbb{Z}}$. Now let $B \subseteq V^{\mathbb{Z}}$ be the set of all such sequences that are bounded. That is, if $x = \{x_n\}_{n \in \mathbb{Z}} \in B$, then there is a $b \in V$ with $\sup\{|x_n|\}_{n \in \mathbb{Z}} \preceq b$. This set is well defined, since by assumption $V$ is Dedekind complete. We set the product topology on $V^{\mathbb{Z}}$, where $(V, T_{ord})$ has the order topology, and we consider $B \subseteq V^{\mathbb{Z}}$ with the subspace topology induced by the topology of $V^{\mathbb{Z}}$. It is worth mentioning that apart from the topology considered in our problem, one can put a norm on $B$, by setting $\|x\|_\infty = \sup_{n \in N} |x_n|$. This norm induces a metric $d(x, y) = \|x - y\|_\infty$ and hence another topology.

**Definition 2.1:** A *Bounded Linear Function (B.L.F.)* on the set $B$ is a map $g : B \to V$ such that for every $x = \{x_n\}_{n \in \mathbb{Z}}$, $y = \{y_n\}_{n \in \mathbb{Z}} \in B$ and $\lambda \in \mathbb{R}$, $g(x + y) = g(x) + g(y)$, $g(\lambda x) = \lambda g(x)$ for and $\inf_{n \in \mathbb{Z}} x_n \preceq g(x) \preceq \sup_{n \in \mathbb{Z}} x_n$.

Notice that the *inf* and the *sup* are well defined, since $\{x_n\}_{n \in \mathbb{Z}}$ is bounded from above and below and $(V, T_{ord})$ is Dedekind complete.

**Example 2.2:** Consider the set $B$ of doubly bounded real valued sequences.
(i) Define $g_1 : B \to \mathbb{R}$ by $g_1(a) = \sum_{n \in \mathbb{Z}} w_n a_n$, where $w_n \in \mathbb{R}$ is a weight such that $\sum_{n \in \mathbb{Z}} w_n = 1$. Then, according to Definition 2.1, $g_1$ is a *B.L.F.* on *B*.

(ii) Define $g_2 : B \to \mathbb{R}$ by $g_2(a) = \sum_{n \in K} a_n$, where $K \subset \mathbb{Z}$ is a finite subset of $\mathbb{Z}$. Then, according to Definition 2.1, $g_2$ is a mean on $B$.

**Remark:** Recall that a permutation of a set G is any bijective function taking G onto G; and the set of all such functions forms a group under function composition, called the symmetric group on G, and written as Sym(G). In group theory Cayley's theorem states that every group G is isomorphic to a subgroup of the symmetric group acting on G (see [5]). This can be understood as an example of the group action of G on the elements of G. The study of the invariance of a function under any Group action of G can be reduced up to isomorphism to the study of the invariance of the Group action of Sym(G) on G.

Let $B \subseteq V^{\mathbb{Z}}$ be the set of doubly infinite bounded sequences in $(V, T_{ord})$. Consider the *symmetric group Sym(B)*. The action of any element $S$ of $Sym(B)$ on $x \in B$ is defined as $S(x_i) = x_{j_i}$, for all $i \in \mathbb{Z}$, where $x = \{x_n\}_{n \in \mathbb{Z}} \in B$, $j_i \in \mathbb{Z}$, and $x_{j_i}$ is the permutation of $x_i$ under the action of $S \in Sym(B)$. To see this, note that if $x = \{x_n\}_{n \in \mathbb{Z}}$, $y = \{y_n\}_{n \in \mathbb{Z}}$ are elements of $B$, then the elements of one sequence are permutations of the elements of the other. That is, for all $i \in \mathbb{Z}$, $y_i = x_{j_i}$ for some $j_i \in \mathbb{Z}$, where $x_{j_i}$ is the permutation of $x_i$ under the action of $S \in Sym(B)$. Hence for any $y \in B$ there exists a unique $x \in B$ such that $S(x) = y$ ($S$ are bijections taking $B$ onto $B$.)

Now note that the two *B.L.F.* in the Example 2.2 are fundamentally different, in the sense that the second, in contrast to the first one, is in addition *permutation invariant* according to the following definition.

**Definition 2.3:** A function $g : B \to V$ is said to be group invariant under $Sym(B)$, written as $Sym(B)$-invariant, if for any $S \in Sym(B)$ and for any $x \in B$ and it holds that $g(S(x)) = g(x)$.

Note also that the second function in Example 2.2, as being $Sym(B)$-invariant, can be considered as an averaging procedure, where the average/mean value is not affected by the ordering of the elements of the sequence.

As we have already mentioned, due to the significance of Cayley's theorem we can reduce a problem of existence of *group invariant* functions to problems regarding existence of functions that are $Sym(B)$-invariant (when groups are considered up to group isomorphisms). The question that we discuss here now becomes concrete and it is the following: Does such a $Sym(B)$-invariant bounded linear function on $B \subseteq V^{\mathbb{Z}}$ exist? And if it does, can we provide a method of constructing it? We pursue this question in the following section and with the aid of *Tychonoff's Theorem*

*(which is equivalent to the Axiom of Choice* [7]) we show that such a function exists on $B \subseteq V^{\mathbb{Z}}$. Since our proof is purely an existential one, as it uses an equivalent statement of the Axiom of Choice, most probably a constructive method is hard to be given.

## 3. The Main Result

Here we present the main result of this article; that is, the existence of *bounded linear functions* (B.L.F.) on $B \subseteq V^{\mathbb{Z}}$ that are invariant under group actions of $Sym(B)$ on $B \subseteq V^{\mathbb{Z}}$. At first we construct "approximately" $Sym(B)$-invariant *B.L.F.* and then we show the existence of $Sym(B)$-invariant *B.L.F.* from the approximately invariant ones.

Next we present some preliminary results necessary for the proof of the main theorem.

**Definition 3.1:** A sequence $\{g_n\}_{n \in \mathbb{N}}$ of *B.L.F.* is approximately $Sym(B)$-invariant if for any $S \in Sym(B)$ and all $x \in B$ $\lim_{n \to \infty} g_n(S(x)) = g_n(x)$.

The existence of such approximate *B.L.F.* is guaranteed by the following Lemma.

**Lemma 3.2:** There is a sequence $\{g_n\}_{n \in \mathbb{N}}$ of *B.L.F.* such that for any $S \in Sym(B)$ and $x \in B$, $\lim_{n \to \infty} g_n(S(x)) = g_n(x)$.

*Proof:* Let $S \in Sym(B)$. For $x = \{x_n\}_{n \in \mathbb{Z}} \in B$, let $g_n(x) = \frac{1}{n^2} \sum_{i=1}^{n} x_i$. Then, by definition 2.3,

$$|g_n(S(x)) - g_n(x)| = \sum_{k=1}^{n} \frac{|x_{j_k} - x_k|}{n^2} \leq \sum_{k=1}^{n} \frac{2|x_k|}{n^2} \leq \frac{2n\|x\|}{n^2} = \frac{2}{n}\|x\|, \text{ where } \|x\| = \sup_{n \in \mathbb{Z}}|x_n| < \infty, \text{ since}$$

$x = \{x_n\}_{n \in \mathbb{Z}}$ is bounded.
Thus, $\lim_{n \to \infty} |g_n(S(x)) - g_n(x)| = 0$. $\square$

How to conclude the existence of $Sym(B)$-invariant B.L.F.? Interestingly, the proof of the existence of $Sym(B)$-invariant functions will be non-constructive. Basically, we consider the set $\mathfrak{M} \subseteq V^B$ of all B.L.F. on $B$ and we topologize it in such a way that $\mathfrak{M}$ becomes a compact space. To do this we use Tychonoff's theorem. Then we show that the sequence $\{g_n\}_{n \in \mathbb{N}}$ of the "*approximately*" $Sym(B)$-invariant B.L.F. has a cluster point $g$, which is $Sym(B)$-invariant. This $g$ is the group invariant B.L.F. that we are looking for.

We view $\mathfrak{M}$ as a subset of $V^B$, by identifying a B.L.F. $g$ with the indexed family $(g(x))_{x\in B}$. We then give $V^B$ the product topology and $\mathfrak{M}$ the subspace topology. Under this topology, the evaluation maps (projections) $g \mapsto g(x)$, $x \in B$, are continuous maps from $\mathfrak{M}$ to $V$ for each $x \in B$. (Actually this is the appropriate topology that makes the evaluation maps continuous, see [8]). Notice also that the definition of the *B.L.F.* actually makes $\mathfrak{M}$ a subset of a product of ordered closed intervals[2]. That is, $\mathfrak{M} \subseteq \prod_{x\in B}[m(x), M(x)]$, where $m(x) = \inf_{n\in\mathbb{Z}} x_n$ and $M(x) = \sup_{n\in\mathbb{Z}} x_n$.

For the proof of the proposition 3.4, we need the following a special case of Tychonoff's theorem.

**Lemma 3.3:** Let $I$ be an indexed set and assume that $(X_i, T_i)$, $i \in I$, are compact Hausdorff spaces. Then $\prod_{i\in I}(X_i, T_{pr})$ with the product topology is also compact.

**Proposition 3.4:** The set $\mathfrak{M} \subseteq V^B$ of all B.L.F on $B$ is a compact space.

*Proof*: We show that $\mathfrak{M}$ is a closed subset of a compact space. Note that $[m(x), M(x)] \subseteq (V, T_{ord})$ is an ordered closed and bounded interval of the totally ordered Dedekind complete space $(V, T_{ord})$, that it has the least upper bound property. Hence, the generalized Heine-Borel theorem implies that $[m(x), M(x)]$ is also compact (see [8] Thm. 27.1). Moreover, $[m(x), M(x)]$ with the subspace topology induced by the $T_{ord}$ on $V$ is a Hausdorff space. Hence, from Lemma 3.3, $X = \prod_{x\in B}[m(x), M(x)]$ is a compact space. We show that $\mathfrak{M} = \bigcap F^{-1}(0)$, for various continuous maps $F: X \to V$. Here recall that the evaluation maps $g \mapsto g(x)$, $x \in B$, are continuous in this topology and the elements $g \in \mathfrak{M}$ are *B.L.F.*. That is, $\forall x, y \in B, \forall \lambda \in \mathbb{R}$,

$$g(x+y) - g(x) - g(y) = 0$$
$$g(\lambda x) - \lambda g(x) = 0$$
(*)

where, for every $x = \{x_n\}_{n\in\mathbb{Z}} \in B$, $\inf_{n\in\mathbb{Z}} x_n \preceq g(x) \preceq \sup_{n\in\mathbb{Z}} x_n$.

This exhibits $\mathfrak{M}$ as an intersection of sets of the form $F^{-1}(\{0\})$. That is, $\mathfrak{M} = \bigcap F^{-1}(\{0\})$ for various continuous maps $F: X \to V$ which are defined via the above continuous evaluation maps $g \mapsto g(x)$ as in (*). Since $F$ is continuous, then $F^{-1}(\{0\})$ is closed (because $\{0\}$ is closed in $(V, T_{ord})$). Thus, as $\mathfrak{M}$ is the intersection of closed sets, it is closed. □

The next lemma is a version of the Bolzano-Weierstrass property for compact spaces.

---

[2] which by Lemma 1.3 are also normed closed.

**Lemma 3.5:** If $X$ *is* a compact space then every sequence $\langle x_n \rangle_{n \in \mathbb{N}}$ in $X$ has a cluster point. That is, there is a point $x \in X$ such that for every open neighborhood $U_x$ of $x \in X$, $U_x \cap \langle x_n \rangle_{n \in \mathbb{N}} \in X$ contains infinitely many points.

**Proof:** For each $x$ that is not a cluster point, there is a neighborhood $U_x$ of $x \in X$ that contains only finitely many points from the sequence $\langle x_n \rangle_{n \in \mathbb{N}}$. Clearly no finite sub-collection of such neighborhoods of $x$ could cover $X$. (Any finite union of such neighborhoods of $x$ contains only finitely many points from $\langle x_n \rangle_{n \in \mathbb{N}} \in X$ and hence cannot cover $\langle x_n \rangle_{n \in \mathbb{N}}$ and thus $X$). Hence, the family $\{U_x\}_{x \in X}$ does not cover $X$ [3]. Since $\bigcup_{x \in X} U_x$ contains all the points of $X$ that are <u>not cluster points,</u> and since $\bigcup_{x \in X} U_x$ does not cover $X$, then $X$ must contain at least one cluster point of $\langle x_n \rangle_{n \in \mathbb{N}}$. Thus, $\langle x_n \rangle_{n \in \mathbb{N}}$ has a cluster point in $X$. □

Now it follows the proof of the main result of this article.

**Theorem 3.6 (Main):** There exists a *B.L.F*, $g \in \mathfrak{M}$, which is *Sym(B)*-invariant; that is, it is invariant under group actions of the symmetric group of $B$ on $B \subseteq V^{\mathbb{Z}}$.

*Proof*: Applying Lemma 3.5 on $\mathfrak{M}$ (which is compact by Proposition.3.4) and on the sequence of "approximately" *Sym(B)*-invariant B.L.F., $\{g_n\}_{n \in \mathbb{N}}$ (see definition 3.1), we conclude that $\{g_n\}_{n \in \mathbb{N}}$ has a cluster point $g$. We show that this $g$ is *Sym(B)*-invariant. Let $S \in Sym(B)$, $x \in B$ arbitrary, and consider the neighborhood $U_g$ of $g$:

$$U_g = \{k \in \mathfrak{M} : \forall x \in B, \; |k(x) - g(x)| < \varepsilon \text{ and } |k(S(x)) - g(S(x))| < \varepsilon \}$$

Actually this is a basic neighborhood defined by imposing conditions on two coordinates of a general element $k \in \mathfrak{M} \subseteq V^B$. Since $g$ is a cluster point of $\{g_n\}_{n \in \mathbb{N}}$, then $g_n \in U_g$ for infinitely many $n$. Let $\varepsilon > 0$ (arbitrary), then $\forall n' \in \mathbb{N} \; \exists n \geq n'$ so that $g_n \in U_g$. We can take $n'$ large enough such that we have (see also Lemma 3.2) $|g_n(S(x)) - g_n(x)| < \varepsilon$, $g_n \in U_g$. Thus,

$$|g(S(x)) - g(x)| \leq |g(S(x)) - g_n(S(x))| + |g_n(S(x)) - g(x)| + |g_n(x) - g(x)| < 3\varepsilon$$

Since $\varepsilon > 0$ is arbitrary, $g(S(x)) = g(x)$.

---

[3] Recall the definition of Compactness.

## 4. Conclusion

In this paper we considered Dedekind-complete totally ordered Riesz Spaces with the order topology, $(V, T_{ord})$. In this setting, we have shown how one can derive the existence of Bounded Linear Functions that are invariant under group actions of the symmetric group of $B$ on $B$, where $B$ is the set of all bounded functions on $(V, T_{ord})$. As we have seen, the main proof of this result relies on the Heine-Borel Theorem. Recall that a topological vector space has the Heine–Borel property, if and only if, every closed and bounded subset of it is compact. Therefore, it is of interest to study the problem of existence of group invariant bounded linear functions whenever we have an ordered structure on general topological vector spaces that have the Heine–Borel property. It is worth mentioning that many metric spaces fail to have the Heine–Borel property. For instance, the metric space of rational numbers (or indeed any incomplete metric space) fails to have the Heine–Borel property. Complete metric spaces may also fail to have the property. For example, no infinite-dimensional Banach space has the Heine–Borel property. On the other hand, some infinite-dimensional Fréchet spaces do have the Heine–Borel property. The space $C^{\infty}(K)$ of smooth functions on a compact set $K \subset \mathbb{R}^n$, considered as a Fréchet space, has the Heine–Borel property, as can be shown by using the Arzelà–Ascoli theorem. This extends to any nuclear Fréchet space (see [3]).

**Conflict of Interests**
The author declares that there is no conflict of interest regarding the publication of this paper.